\documentclass[12pt,a4paper]{article}
\usepackage{pifont}
\usepackage{amsfonts}
\usepackage{mathrsfs}
\usepackage{amssymb}
\usepackage{amsmath}

\newtheorem{theorem}{Theorem}[section]
\newtheorem{definition}[theorem]{Definition}
\newtheorem{lemma}[theorem]{Lemma}
\newtheorem{corollary}[theorem]{Corollary}

\newtheorem{example}[theorem]{Example}

\makeatletter

\newcommand{\Rmnum}[1]{\expandafter\@slowromancap\romannumeral #1@}
\makeatother

\begin{document}

\title{Gr\"{o}bner-Shirshov bases for semirings\footnote{Supported by the NNSF of China (No. 11171118),
the Research Fund for the Doctoral Program of Higher Education of
China (No. 20114407110007) and the NSF of Guangdong Province (No.
S2011010003374).}}
\author{
L. A. Bokut\footnote {Supported by RFBR 09-01-00157,
LSS--3669.2010.1 and SB RAS Integration grant No. 2009.97 (Russia)
and Federal Target Grant ¡°Scientific and educational personnel of
innovation Russia¡± for 2009-2013
(government contract No. 02.740.11.5191).} \\
{\small \ School of Mathematical Sciences, South China Normal
University}\\
{\small Guangzhou 510631, P. R. China}\\
{\small Sobolev Institute of Mathematics, Russian Academy of
Sciences}\\
{\small Siberian Branch, Novosibirsk 630090, Russia}\\
{\small Email: bokut@math.nsc.ru}\\
\\
 Yuqun
Chen\footnote {Corresponding author.} \  and Qiuhui Mo\\
{\small \ School of Mathematical Sciences, South China Normal
University}\\
{\small Guangzhou 510631, P. R. China}\\
{\small Email: yqchen@scnu.edu.cn}\\
{\small scnuhuashimomo@126.com}}

\date{{\bf In Memorial of Jean-Louis Loday}}

\maketitle \noindent\textbf{Abstract:} In the paper, we establish
Gr\"obner-Shirshov bases for semirings and commutative semirings. As
applications, we obtain Gr\"obner-Shirshov bases and A. Blass's
(1995) and M. Fiore -T. Leinster's (2004) normal forms of the
semirings $\mathbb{N}[x]/(x=1+x+x^2)$ and $\mathbb{N}[x]/(x=1+x^2)$
with one generator $x$ and one defining relation, correspondingly.

\noindent \textbf{Key words:}  Gr\"{o}bner-Shirshov basis, semiring,
 congruence, normal form.

\noindent \textbf{AMS 2010 Subject Classification}: 16Y60, 16S15,
13P10

\section{Introduction}
Gr\"{o}bner bases and Gr\"{o}bner-Shirshov bases  were invented
independently by A.I. Shirshov for ideals of free (commutative,
anti-commutative) non-associative algebras \cite{Sh62a,Shir3}, free
Lie algebras \cite{Sh62b,Shir3} and implicitly free associative
algebras \cite{Sh62b,Shir3}  (see also \cite{Be78,Bo76}), by H.
Hironaka \cite{Hi64} for ideals of the power series algebras (both
formal and convergent), and by B. Buchberger \cite{Bu70} for ideals
of the polynomial algebras.

Gr\"{o}bner bases and Gr\"{o}bner-Shirshov bases theories have been
proved to be very useful in different branches of mathematics,
including commutative algebra and combinatorial algebra, see, for
example, the books
 \cite{AL, BKu94, BuCL, BuW, CLO, Ei}, the papers \cite{Be78, Bo72,Bo76},
 and the surveys \cite{BC,BCS, BFKK00, BK03, BK05}.

Up to now, different versions of Composition-Diamond lemma are known
for the following classes of algebras apart those mentioned above:
(color) Lie super-algebras \cite{Mik89, Mik92,Mik96}, Lie
$p$-algebras \cite{Mik92}, associative conformal algebras
\cite{BFK}, modules \cite{Chi,KL} (see also \cite{CCZ}),
right-symmetric algebras \cite{BCLi08}, dialgebras \cite{BCL08},
associative algebras with multiple operators \cite{BCQ08},
Rota-Baxter algebras \cite{BCD08}, and so on.

It is well known Shirshov's result \cite{Shir1, Shir3} that every
finitely or countably generated Lie algebra over a field $k$ can be
embedded into a two-generated Lie algebra over $k$.
 Actually, from the technical
point of view, it was a beginning of the Gr\"{o}bner-Shirshov bases
theory for Lie algebras (and associative algebras as well). Another
proof of the result using explicitly Gr\"{o}bner-Shirshov bases
theory is refereed
 to  L.A. Bokut, Yuqun Chen  and
Qiuhui Mo \cite{BCM}.

A.A. Mikhalev and  A.A. Zolotykh \cite{MZ} prove the
Composition-Diamond lemma for a tensor product of a free algebra and
a polynomial algebra, i.e. they establish Gr\"{o}bner-Shirshov bases
theory for
 associative algebras over
 a commutative algebra. L.A. Bokut, Yuqun Chen and Yongshan Chen \cite{BCC08} prove the
Composition-Diamond lemma for a tensor product of two free algebras.
Yuqun Chen, Jing Li and Mingjun Zeng \cite{CJZ} prove the
Composition-Diamond lemma for a tensor product of a non-associative
algebra and a polynomial algebra.

L.A. Bokut, Yuqun Chen and Yongshan Chen \cite{BCC11} establish the
Composition-Diamond lemma for Lie algebras over a polynomial
algebra, i.e. for ``double free" Lie algebras. It provides a
Gr\"{o}bner-Shirshov bases theory for Lie algebras over a
commutative algebra. Yuqun Chen and Yongshan Chen \cite{CC12}
establish the Composition-Diamond lemma for matabelian Lie algebras.

In this paper, we establish Gr\"{o}bner-Shirshov bases for semirings
and  commutative semirings. We show that for a given monomial
ordering on the free (commutative) semiring, each ideal of the free
(commutative) semiring algebra has the unique reduced
Gr\"{o}bner-Shirshov basis.

In 2004, M. Fiore and T. Leinster \cite{mar} find a strongly
normalizing reduction system and a normal form of the semiring
$\mathbb{N}[x]/(x=1+x+x^2)$ where $\mathbb{N}$ is the set of natural
numbers which is regarded as a semiring and $(x=1+x+x^2) \
((x=1+x^2)) $ is the congruence on the semiring $\mathbb{N}[x]$
generated by $x=1+x+x^2 \ (x=1+x^2) $.  In 1995, A. Blass
\cite{blass} finds a normal form of the semiring
$\mathbb{N}[x]/(x=1+x^2)$. Now, we use the Composition-Diamond lemma
for commutative semirings to find  Gr\"{o}bner-Shirshov bases  and
Fiore-Leinster's and Blass's normal forms for the above semirings
respectively. Also we show that each congruence of the semiring
$\mathbb{N}$ is generated by one element and that the commutative
semiring $\mathbb{N}[x]$ is not Noetherian.

\section{Free semiring}

Let $A$ be an $\Omega$-system,  $
\Omega=\bigcup_{n=1}^{\infty}\Omega_{n}$ where $\Omega_{n}$ is the
set of $n$-ary operations, for example, ary $(\delta)=n$ if
$\delta\in \Omega_n$.  We will call $A$ an $\Omega$-groupoid and
each $\omega$ would be called a product. Let $k$ be a field and $kA$
the groupoid algebra over $k$, i.e. the $k$-linear space  with a
$k$-basis $A$ and linear multiple products $\omega\in \Omega$ that
are extended  by linearity from $A$ to $kA$. Such $kA$ is called an
$\Omega$-algebra. If $\rho$ is a congruence relation on $A$
generated by pairs $(a_i,b_i)$, $i\in I$, then, as
$\Omega$-algebras,
$$
k(A/\rho)\cong kA/Id(a_i-b_i, i\in I),
$$
where $Id(S)$ is the $\Omega$-ideal of $kA$ generated by $S$. It
means that any ``monomial" linear basis of $kA/Id(a_i-b_i, i\in I)$,
i.e. a basis that consists of elements of $A$, is a set of normal
forms for $A$.

Now let $(A, \circ, \cdot, \theta, 1)$ be a semiring, i.e. $(A,
\circ, \theta)$ is a commutative monoid, $(A, \cdot, 1)$ is a
monoid,  $\theta\cdot a=a\cdot \theta=\theta$ for any $a\in A$, and
$\cdot$ is distributive relative to $\circ$ from left and right.

Some people call ``rig" instead of ``semiring".

The class of  semirings is a variety.  So a free semiring
$Rig\langle X\rangle$ generated by a set $X$ is defined as usual as
for any variety of universal systems. Let $(X^*,\cdot,1)$ be the
free monoid generated by $X$. If one fixes some linear ordering $<$
on  $X^*$, then any element of  $Rig\langle X\rangle$ has a unique
form  $w=u_1\circ u_2\circ\cdots\circ u_n$, where $u_i\in X^*,
u_1\leq u_2\leq \ldots \leq u_n,\ n\geq0$ and $w=\theta$ if $n=0$.

For any $ u_1, u_2,  \ldots, u_n\in X^*$ we have $u_1\circ u_2 \circ
 \cdots \circ u_n\in Rig\langle X\rangle$ and $u_1\circ u_2
 \circ \cdots \circ u_n=u_{i_1}\circ u_{i_2}
\circ \cdots \circ u_{i_n}$ where $\{u_1, u_2,  \ldots,
u_n\}=\{u_{i_1}, u_{i_2}, \ldots, u_{i_n}\}$ and $u_{i_1}\leq
u_{i_2}\leq \ldots\leq u_{i_n}$.

$Rig\langle X\rangle=(X^*,\cdot,\circ)$ is a $\Omega$-groupoid,
where $\Omega=\{\cdot,\circ\}$. From now on, we assume that
$\Omega=\{\cdot,\circ\}$.

\ \

Let $k$ be a field. We call the groupoid algebra $kRig\langle
 X\rangle$ a semiring algebra.  Any
element in $Rig\langle X\rangle$ is called a monomial and any
element in $kRig\langle
 X\rangle$ is called a polynomial.

For any $u_1\circ u_2  \circ \cdots \circ u_n\in Rig\langle
X\rangle$, $u_i\in
 X^*$, we denote $|u|_{\circ}=n$.

For any $u\circ u \circ u\circ \cdots \circ u$ where $|u\circ u
\circ u\circ \cdots \circ u|_{\circ}=n$, $u\in
 X^*$, we will denote it by $u^{\circ n}$.

For any  $u=w_1\circ w_2\circ \cdots \circ w_m \circ u_{m+1} \circ
\cdots \circ u_n$, $v=w_1\circ w_2\circ \cdots \circ w_{m}\circ
v_{m+1}\circ\cdots \circ v_t$, where $w_l,u_i, v_j\in X^*$,  such
that
$$
u_i\neq v_j \ \  for \ \ any  \ \ i=m+1, \ldots,n, j=m+1,
\ldots,t,
$$
we denote
$$
lcm_\circ (u, v)=w_1\circ w_2\circ \cdots \circ w_m \circ u_{m+1}
\circ \cdots \circ u_n\circ v_{m+1}\circ\cdots \circ v_t
$$
the least common multiple of $u$ and $v$ in $Rig\langle X\rangle$
with respect to $\circ$.

Throughout this paper, we denote $\mathbb{N}$ the set of natural
numbers, $Rig\langle X|S\rangle$ the semiring with generators $X$
and relations $S$.

\ \

We want to create Gr\"{o}bner-Shirshov bases theory for $Rig\langle
X\rangle$. As for semigroups or groups, it is enough to create
Gr\"{o}bner-Shirshov bases theory for the algebra $kRig\langle
X\rangle$.  Let us remind that $\theta$ and 0 are different elements
of $kRig\langle X\rangle$ since $\theta\in Rig\langle X\rangle$ and
$0\notin Rig\langle X\rangle$.

\section{Composition-Diamond lemma for semirings}

Let  $\star\notin X$. By a $\star$-monomial we mean a monomial in
$Rig\langle
 X\cup \star\rangle$ with only one occurrence of $\star$. Let
u be a  $\star$-monomial  and $s\in kRig\langle
 X\rangle$. Then we
call
$$
u|_{s}=u|_{\star\mapsto s}
$$
an $s$-monomial.  For example, if
$$
u=x_1x_2\circ x_2\star x_3\in Rig\langle
 X\cup \star\rangle
$$
then
$$
u|_s=u|_{\star\mapsto s}=x_1x_2\circ x_2s x_3.
$$

Let $>$ be any monomial ordering on $Rig\langle
 X\rangle$, i.e. $>$ is a well ordering  such that for any $v, w\in
Rig\langle X\rangle$ and $u$ a $\star$-monomial,
$$
w>v\Rightarrow u|_w>u|_v.
$$

For every polynomial $f\in kRig\langle
 X\rangle$, $f$
has the leading monomial $\bar{f}$. If the coefficient of $\bar{f}$
is $1$, then we call $f$ to be monic.

For any set $S\subseteq kRig\langle
 X\rangle $, we say $S$  monic if any $s\in S$ is monic.

\begin{definition}\label{d1}
Let $>$ be a monomial ordering on $Rig\langle
 X\rangle$. Let $f, g$ be two monic polynomials in $kRig\langle
 X\rangle$.
\begin{enumerate}
\item[(I)]If there exist $a,b\in
 X^*$, such that $|lcm_\circ(\bar{f}a, b\bar{g}) |_{\circ}<|\bar{f}a|_{\circ}+ |b\bar{g}|_{\circ}$
then  we call $(f,g)_{w}=fa\circ u-bg\circ v$ the intersection
composition of $f$ and $g$ with respect to $w$ where
$w=lcm_\circ(\bar{f}a, b\bar{g})=\bar{f}a\circ u=b\overline{g}\circ
v$.
\item[(II)] If there exist $a,b\in
 X^*$, such that $|lcm_\circ(\bar{f}, a\bar{g}b)|_{\circ}<|\bar{f}|_{\circ}+ |a\bar{g}b|_{\circ}$
then  we call $(f,g)_{w}=f\circ u-agb\circ v$ the inclusion
composition of  $f$ and $g$ with respect to $w$ where
$w=lcm_\circ(\bar{f}, a\bar{g}b)=\bar{f}\circ u=a\overline{g}b\circ
v$.
\end{enumerate}
\end{definition}

In the above definition, $w$ is called the ambiguity of the
composition. Clearly,
$$
(f,g)_w\in Id(f,g) \ \ \ \mbox{ and }  \ \ \ \overline{(f,g)_w}< w,
$$
where $Id(f,g)$ is the ideal of $kRig\langle
 X\rangle$ generated
by $f,\ g$.

\ \

\noindent{\bf Remark}: We regard $kRig\langle
 X\rangle$ as an $\Omega$-algebra. In the Definition \ref{d1}, the
 ideal  $Id(f,g)$ means $\Omega$-ideal. In this paper, the ideal of $kRig\langle
 X\rangle$ will be $\Omega$-ideal.

\ \

Let $f, g $ be polynomials and $g$ monic with
$\bar{f}=a\overline{g}b\circ u$ for some $a,b\in
 X^*, u\in Rig\langle X\rangle $.  Then
the transformation
$$
f\rightarrow f-\alpha a\overline{g}b\circ u
$$
is called the elimination of the leading term (ELT) of $f$ by $g$,
where $\alpha$ is the coefficient of the leading term of $f$.

\begin{definition}
Suppose that $w$ is a monomial, $S$  a set of monic polynomials in
$kRig\langle  X\rangle$ and $h$  a polynomial. Then  $h$ is trivial
modulo $(S,w)$, denoted by $h\equiv0 \ mod(S,w)$, if
$h=\sum\limits_{i}\alpha_ia_is_ib_i\circ u_i$, where each
$\alpha_{i}\in{k}$, $a_i, b_i\in  X^*, u_i\in Rig\langle X\rangle$,
$s_i\in{S}$ and $a_i\overline{s_i}b_i\circ u_i<w$.

The set $S$ is called a Gr\"{o}bner-Shirshov basis in $kRig\langle
 X\rangle$ if
any composition in $S$ is trivial modulo $S$ and corresponding to
$w$.
\end{definition}

A set $S$ is called a minimal Gr\"{o}bner-Shirshov basis in
$kRig\langle
 X\rangle$ if $S$ is  a Gr\"{o}bner-Shirshov basis in $kRig\langle
 X\rangle$ and for any $f,g\in S$ with $f\neq g$, $\not\exists \ a, b\in  X^*, u\in Rig\langle X\rangle, s.t., \  \overline{f}=a\overline{g}b\circ
 u$.

 Denote
$$
Irr(S) = \{ w\in  Rig\langle X\rangle \ | \ w \neq
{a\overline{s}b\circ u}
 \mbox{ for  any} \ a,b\in
 X^*, u\in Rig\langle X\rangle, s\in S\}.
$$

A  Gr\"{o}bner-Shirshov basis $S$ in $kRig\langle
 X\rangle$ is reduced if for any $s\in S, \ supp(s)\subseteq
 Irr(S-\{s\})$, where $supp(s)=\{u_1,u_2,\dots,u_n\}$ if $s=\sum_{i=1}^n\alpha_iu_i,\ 0\neq\alpha_i\in k,\ u_i\in Rig\langle
 X\rangle$.

If the set $S$ is  a  Gr\"{o}bner-Shirshov basis in $kRig\langle
 X\rangle$, then we call also $S$ is  a  Gr\"{o}bner-Shirshov basis
 for the ideal $Id(S)$ or the algebra $kRig\langle
 X|S\rangle:=kRig\langle
 X\rangle/Id(S)$.

Let $I$ be an ideal of $kRig\langle X\rangle$. Then there exists
uniquely the reduced Gr\"{o}bner-Shirshov basis $S$ for $I$, see
Theorem \ref{t5}.

 \ \

If a subset $S$ of $kRig\langle
 X\rangle$ is not a Gr\"{o}bner-Shirshov basis for $Id(S)$
then one can add to $S$ a nontrivial composition $(f, g)_w$ of $f,
g\in S$ and continue this process repeatedly (actually using the
transfinite induction) in order to obtain a set $S^{comp}$ of
generators of $Id(S)$ such that $S^{comp}$ is a Gr\"{o}bner-Shirshov
basis in  $kRig\langle
 X\rangle$. Such a process is called Shirshov algorithm.

 \ \

Suppose that $S=\{u_i-v_i\ |\ i\in I\}$ where for any $i$, $u_i,
v_i\in Rig\langle
 X\rangle$. In this case we call $u_i-v_i$ a semiring relation.  It is clear that
 if $S$ is a set of  semiring relations then so is $S^{comp}$. In
 order to find a normal form of the semiring $Rig\langle
 X|S\rangle=Rig\langle
 X\rangle/\rho(S)$, where $\rho(S)$ is the congruence of $Rig\langle
 X\rangle$ generated by the set $\{(u_i,v_i)|i\in I\}$, it is enough to find a
 monomial $k$-basis of the semiring algebra $kRig\langle
 X|S\rangle$.

\begin{lemma}\label{l1}
Let $S$ be a Gr\"{o}bner-Shirshov  basis  in  $kRig\langle
 X\rangle$
 and $s_1, s_2\in S$. If
 $w={a\overline{s_1}b\circ u}={c\overline{s_2}d\circ v}$ for some $a,b,c,d\in
 X^*, u,v\in Rig\langle X\rangle$, then
$$
as_1b\circ u\equiv cs_2d\circ v\ \ mod (S,w).
$$
\end{lemma}

{\bf Proof:}  There are two cases to consider.

(I)\ \ $lcm_\circ(a\overline{s_1}b,
c\overline{s_2}d)=a\overline{s_1}b \circ c\overline{s_2}d$ which
means there exists $u_1\in Rig\langle X\rangle$ such that
$u=u_1\circ c\overline{s_2}d,\ v=a\overline{s_1}b\circ u_1$.

Then
\begin{eqnarray*}
&& as_1b\circ u- cs_2d\circ v \\
&=& as_1b\circ u_1\circ
c\overline{s_2}d- cs_2d\circ a\overline{s_1}b\circ u_1 \\
&=&as_1b\circ c\overline{s_2}d\circ u_1-as_1b\circ cs_2d\circ u_1+as_1b\circ cs_2d\circ u_1-a\overline{s_1}b\circ c s_2d\circ  u_1\\
&=&-as_1b\circ c(s_2-\overline{s_2})d\circ
u_1+a(s_1-\overline{s_1})b\circ cs_2d\circ u_1.
\end{eqnarray*}

Since $\overline{s_2-\overline{s_2}}<\overline{s_2}$ and
$\overline{s_1-\overline{s_1}}<\overline{s_1}$, we have
$$
\overline{as_1b\circ c(s_2-\overline{s_2})d\circ u_1}<
a\overline{s_1}b\circ c\overline{s_2}d\circ u_1 =w
$$
and
$$
\overline{a(s_1-\overline{s_1})b\circ cs_2d\circ
u_1}<a\overline{s_1}b\circ c\overline{s_2}d\circ u_1=w.
$$

It follows that

$$
as_1b\circ u\equiv cs_2d\circ v\ mod (S,w).
$$

(II) \ \ $|lcm_\circ(a\overline{s_1}b, c\overline{s_2}d)|_\circ
<|a\overline{s_1}b|_\circ + |c\overline{s_2}d|_\circ$, i.e.
$\overline{s_1}=u_1\circ u_2\circ \cdots \circ u_m \circ u_{m+1}
\circ \cdots \circ u_n$, $\overline{s_2}=v_1\circ v_2\circ \cdots
\circ v_{m}\circ v_{m+1}\circ\cdots \circ v_t$ such that
$$
au_1b=cv_1d,\  au_{2}b=cv_2d,  \ldots,  au_mb= cv_{m}d
$$
and
$$ au_ib\neq cv_jd \ \  for \ \ any  \ \ i=m+1, \ldots,n, \ j=m+1,
\ldots,t,
$$
where $u_i, v_j\in X^*$. In this case, there exists $u'\in
Rig\langle X\rangle$  such that
$$
u=cv_{m+1}d\circ cv_{m+2}d\circ \cdots \circ cv_{t}d\circ u', \ \
v=au_{m+1}b\circ au_{m+2}b\circ\cdots \circ au_{n}b\circ u'.
$$

There are four subcases to consider.

1) $u_1b_1=c_1v_1$ for some $b_1, c_1\in X^*$. In this subcase,
$c=ac_1, b=b_1d$. Then
\begin{eqnarray*}
&& as_1b\circ u- cs_2d\circ v \\
&=&as_1b_1d\circ ac_1v_{m+1}d\circ
ac_1v_{m+2}d\circ \cdots \circ ac_1v_td\circ u' \\
&&-ac_1s_2d\circ au_{m+1}b_1d\circ au_{m+2}b_1d\circ \cdots \circ au_{n}b_1d\circ u'\\
&=&a(s_1b_1\circ c_1v_{m+1}\circ c_1v_{m+2}\circ\cdots\circ
c_1v_t-c_1s_2\circ u_{m+1}b_{1}\circ u_{m+2}b_1\circ \cdots\circ
u_{n}b_1)d\circ u'\\
&=&a( (s_1,s_2)_{w'} )d\circ u'
\end{eqnarray*}
where $w'=\overline{s_1}b_1\circ c_1v_{m+1}\circ
c_1v_{m+2}\circ\cdots\circ c_1v_t=c_1\overline{s_2}\circ
u_{m+1}b_{1}\circ u_{m+2}b_1\circ \cdots\circ u_{n}b_1$.

Since  $S $ is a Gr\"{o}bner-Shirshov basis in  $kRig\langle
 X\rangle$, we have
$$
(s_1,s_2)_{w'}\equiv 0 \  mod (S,w').
$$

Then
\begin{eqnarray*}
&& as_1b\circ u- cs_2d\circ v \\
&=&a((s_1,s_2)_{w'} )d\circ u'\\
&\equiv& 0 \ mod (S,aw'd\circ u') \\
&\equiv& 0 \ mod (S,w).
\end{eqnarray*}

2) $a_1u_1=v_1d_1$ for some $a_1, d_1\in X^*$. This subcase is
similar to subcase 1). We omit the proof.

3) $u_1=c_1v_1d_1$ for some  $c_1, d_1\in X^*$. Then $c=ac_1,
d=d_1b$ and
\begin{eqnarray*}
&& as_1b\circ u- cs_2d\circ v \\
&=& as_1b\circ ac_1v_{m+1}d_1b\circ ac_1v_{m+2}d_1b\circ \cdots \circ  ac_1v_td_1b\circ u' \\
&&-ac_1s_2d_1b\circ au_{m+1}b\circ au_{m+2}b\circ \cdots \circ  au_{n}b\circ u' \\
&=& a(s_1\circ c_1v_{m+1}d_1\circ c_1v_{m+2}d_1\circ \cdots \circ
c_1v_td_1-c_1s_2d_1\circ u_{m+1}\circ u_{m+2}\circ \cdots \circ
u_{n})b\circ u'\\
&=&a( (s_1,s_2)_{w'} )b\circ u'
\end{eqnarray*}
where  $w'=\overline{s_1}\circ c_1v_{m+1}d_1\circ c_1v_{m+2}d_1\circ
\cdots \circ c_1v_td_1=c_1\overline{s_2}d_1\circ u_{m+1}\circ
u_{m+2}\circ \cdots \circ u_{n}$.

Since  $S $ is a Gr\"{o}bner-Shirshov basis in  $kRig\langle
 X\rangle$, we have
$$
(s_1,s_2)_{w'}\equiv 0 \  mod (S,w').
$$
Similar to the subcase 1), we have
\begin{eqnarray*}
&& as_1b\circ u- cs_2d\circ v  \equiv 0 \ mod (S,w).
\end{eqnarray*}

4) $a_1u_1b_1=v_1$ for some  $a_1, b_1\in X^*$. This subcase is
similar to the subcase 3). We omit the proof.

The lemma is proved. \hfill $\blacksquare$

\ \

\begin{theorem}\label{t1}{\em(Composition-Diamond lemma for semirings)}\ \  Let $S$ be a set of monic
polynomials in $kRig\langle
 X\rangle$, $>$ a monomial ordering on $Rig\langle
 X\rangle$  and $Id(S)$ the $\Omega$-ideal of $kRig\langle
 X\rangle$ generated by $S$.  Then the following statements are
equivalent.
 \begin{enumerate}
\item[(1)] $S $ is a Gr\"{o}bner-Shirshov basis in $kRig\langle
 X\rangle$.
\item[(2)] $ f\in Id(S)\Rightarrow \bar{f}=a\overline{s}b\circ u$
for some $a,b\in
 X^*, u\in Rig\langle X\rangle$  and $s\in S$.
\item[(2$^{'}$)]  $f\in Id(S)\Rightarrow
f=\alpha_1a_1s_1b_1\circ u_1+\alpha_2a_2s_2b_2\circ
u_2+\ldots+\alpha_na_ns_nb_n\circ u_n $, where
$a_1\overline{s_1}b_1\circ u_1>a_2\overline{s_2}b_2\circ
u_2>\ldots>a_n\overline{s_n}b_n\circ u_n$, $0\neq\alpha_{i}\in{k},\
a_i, b_i\in X^*,\ u_i\in Rig\langle X\rangle$, $s_i\in{S}$.
\item[(3)] $Irr(S) = \{ w\in  Rig\langle X\rangle \ | \ w \neq {a\overline{s}b\circ u}
 \mbox{ for  any} \ a,b\in
 X^*, u\in Rig\langle X\rangle, s\in S\}$
is a $k$-basis of $kRig\langle
 X|S\rangle=kRig\langle
 X\rangle/Id(S)$.
\end{enumerate}
\end{theorem}
{\bf Proof:} (1)$\Longrightarrow$ (2)\ \ Let  $0\neq f\in Id(S)$.
Then

$$
f=\sum\limits_{i=1}^{n}\alpha_i a_is_ib_i\circ u_i
$$
where each $\alpha_i\in k, a_i, b_i\in X^*, u_i\in Rig\langle
X\rangle$, $s_i\in{S}$.

Let $w_i=a_i\overline{s_i}b_i\circ u_i$ and we arrange this leading
terms in non-increasing ordering by
$$
w_1= w_2=\ldots=w_m >w_{m+1}\geq \ldots\geq w_n.
$$
Now we prove the result by  induction  on $m$.

If $m=1$, then $\bar{f}=a_1\overline{s_1}b_1\circ u_1$.

Now we assume that $m\geq 2$. Then
$$
a_1\overline{s_1}b_1\circ u_1=w_1=w_2=a_2\overline{s_2}b_2\circ u_2.
$$

Since $S $ is a Gr\"{o}bner-Shirshov basis  in $kRig\langle
 X\rangle$,  by Lemma
\ref{l1}, we have
$$
a_2s_2b_2\circ u_2-a_1s_1b_1\circ u_1=\sum\beta_jc_js_jd_j\circ v_j
$$
where each $\beta_j\in k, c_j, d_j\in X^*, v_j\in Rig\langle
X\rangle$, $s_j\in{S}$, and $c_j\overline{s_j}d_j\circ v_j<w_1$.
Therefore, since
$$
\alpha_1a_1s_1b_1\circ u_1+\alpha_2a_2s_2b_2\circ
u_2=(\alpha_1+\alpha_2)a_1s_1b_1\circ u_1+\alpha_2(a_2s_2b_2\circ
u_2-a_1s_1b_1\circ u_1),
$$
we have
$$
f=(\alpha_1+\alpha_2)a_1s_1b_1\circ
u_1+\alpha_2\sum\beta_jc_js_jd_j\circ v_j+
\sum\limits_{i=3}^{n}\alpha_ia_is_ib_i\circ u_i.
$$

If either $m>2$ or $\alpha_1+\alpha_2\neq 0$, then the result
follows from the induction on $m$. If $m=2$ and
$\alpha_1+\alpha_2=0$, then  the result follows from the induction
on $w_1$.

(2)$\Leftrightarrow (2')$ is clear.

(2)$\Longrightarrow$ (3) For any $f\in kRig\langle
 X\rangle$, by the ELTs, we can obtain
 that $f+Id(S)$ can be expressed as a linear combination of  elements of $Irr(S)$.
 Now suppose $\alpha_1u_1+\alpha_2u_2+\ldots+\alpha_nu_n=0$ in
$kRig\langle
 X|S\rangle$ with $u_i\in Irr(S)$,
$u_1>u_2>\ldots>u_n$ and each $\alpha_i\neq 0$.  Then, in
$kRig\langle
 X\rangle$,
$$
g=\alpha_1u_1+\alpha_2u_2+\ldots+\alpha_nu_n\in Id(S).
$$
By (2), we have $u_1=\bar{g}\notin Irr(S)$, a contradiction. So
$Irr(S)$ is $k$-linearly independent. This shows that $Irr(S)$ is a
$k$-basis of $kRig\langle
 X|S\rangle$.

(3)$\Longrightarrow $(2) Let $0\neq f\in Id(S)$. Suppose that
$\bar{f}\in Irr(S)$. Then
$$
f+Id(S)=\alpha (\bar{f}+Id(S))+ \sum \alpha_i(u_i+Id(S)),
$$
where $\alpha,\alpha_i\in k,\ u_i\in Irr(S)$ and $\bar{f}>u_i$.
Therefore, $f+Id(S)\neq 0$, a contradiction. So
$\bar{f}=a\overline{s}b\circ u$ for some $a, b\in X^*, u\in
Rig\langle X\rangle$, $s\in{S}$.

 (2)$\Longrightarrow $(1) \ By the definition of the   composition,
we have  $(f,g)_w\in Id(S)$. If $(f,g)_w\neq 0$, then by (2), $
\overline{(f,g)_w}=a_1\overline{s_1}b_1\circ u_1$ for some $a_1,
b_1\in X^*, u_1\in Rig\langle X\rangle$, $s_1\in{S}$. Let
$$
h=(f,g)_w-\alpha_1a_1{s_1}b_1\circ u_1,
$$
where $\alpha_1$ is the coefficient of   $ \overline{(f,g)_w}$. Then
$\bar{h}< \overline{(f,g)_w}$ and $h\in Id(S)$. By induction, we can
get the result.
 \hfill $\blacksquare$

\ \

Suppose that $>$ is a monomial ordering on $Rig\langle
 X\rangle$ and $I$ an ideal of $kRig\langle
 X\rangle$. Then there exists a Gr\"{o}bner-Shirshov basis $S\subset kRig\langle
 X\rangle$  for the ideal $I=Id(S)$, for example, we may take $S=I$. By Theorem \ref{t1}, we may assume that the leading
 terms of the elements of $S$ are different with each other. For any
 $g\in S$, denote
$$
\Delta_g=\{f\in S|f\neq g \mbox{ and } \overline{f}=a\bar{g}b\circ
u \mbox{ for some } a,b\in
 X^*, u\in Rig\langle X\rangle\}
$$
and $S_1=S-\cup_{g\in S}\Delta_g$.

For any $f\in Id(S)$ we show that there exists an $s_1\in S_1$ such
that $\overline{f}=a\overline{s_1}b\circ u \mbox{ for some } a,b\in
 X^*, u\in Rig\langle X\rangle$.

In fact, by Theorem \ref{t1}, $\overline{f}=a'\bar{h}b'\circ u'
\mbox{ for some } a',b'\in
 X^*, u'\in Rig\langle X\rangle$ and $h\in S$. Suppose that $h\in S-S_1$. Then we have $h\in
\cup_{g\in S}\Delta_g$, say, $h\in \Delta_g$, i.e. $h\neq g \mbox{
and } \overline{h}=a\bar{g}b\circ u \mbox{ for some } a,b\in
 X^*, u\in Rig\langle X\rangle$. We claim that $\bar{h}>\bar{g}$.
 Otherwise, $\bar{h}<\bar{g}$. It follows that   $\bar{h}=a\bar{g}b\circ u>a\bar{h}b\circ
 u$ and so we have an infinite descending chain
$$
\bar{h}>a\bar{h}b\circ u >a^2\bar{h}b^2\circ aub\circ u
>a^3\bar{h}b^3\circ a^2ub^2\circ aub\circ u>\dots
$$
which contradicts that $>$ is well ordered.

Suppose that $g\not\in S_1$. Then
 by the above proof, there exists a $g_1\in S$ such that $g\in
 \Delta_{g_1}$ and $\overline{g}>\overline{g_1}$. Since  $>$ is well
 ordered, there must exist an $s_1\in S_1$ such that
$\overline{f}=a_1\overline{s_1}b_1\circ u_1 \mbox{ for some }
a_1,b_1\in
 X^*, u_1\in Rig\langle X\rangle$.

 Let $f_1=f-\alpha_1a_1s_1b_1\circ u_1$, where $\alpha_1$ is the
 coefficient of the leading term of $f$. Then $f_1\in Id(S)$ and
 $\overline{f}>\overline{f_1}$.

By induction on $\overline{f}$, we know that $f\in Id(S_1)$ and
hence $I=Id(S_1)$. Moreover, by Theorem \ref{t1}, $S_1$ is clearly a
minimal Gr\"{o}bner-Shirshov basis for the ideal $Id(S)$.

Assume that $S$ is  a minimal Gr\"{o}bner-Shirshov basis for the
ideal $I$.

For any $s\in S$, we have $s=s'+s''$, where $supp(s')\subseteq
Irr(S-\{s\}),\ s''\in Id(S-\{s\})$. Since $S$ is  a minimal
Gr\"{o}bner-Shirshov basis, we have $\overline{s}=\overline{s'}$ for
any $s\in S$.

Then $S_2=\{s'|s\in S\}$ is the reduced Gr\"{o}bner-Shirshov basis
for the ideal $I$. In fact, it is clear that $S_2\subseteq Id(S)=I$.
For any $f\in Id(S)$, by Theorem \ref{t1},
$\overline{f}=a_1\overline{s_1}b_1\circ u_1
=a_1\overline{s_1'}b_1\circ u_1\mbox{ for some } a_1,b_1\in  X^*,
u_1\in Rig\langle X\rangle$.

Suppose that $S,\ R$ are two reduced Gr\"{o}bner-Shirshov bases for
the ideal $I$. For any $s\in S$, by Theorem \ref{t1},
$$
\overline{s}=a\overline{r}b\circ u, \
\overline{r}=c\overline{s_1}d\circ v
$$
for some $ a,b,c,d\in
 X^*, u,v\in Rig\langle X\rangle$ and hence $\overline{s}=ac\overline{s_1}db\circ avb\circ
 u$. Since $\bar s\in supp(s)\subseteq Irr(S-\{s\})$, we have
 $s=s_1$. It follows that $a=b=c=d=1$ and $u=v=\theta$ and so
 $\overline{s}=\overline{r}$.

If $s\neq r$ then $0\neq s-r\in I=Id(S)=Id(R)$. By  Theorem
\ref{t1}, $\overline{s-r}=a_1\overline{r_1}b_1\circ
u_1=c_1\overline{s_2}d_1\circ v_1$ for some $ a_1,b_1,c_1,d_1\in
 X^*, u_1,v_1\in Rig\langle X\rangle$ with
 $\overline{r_1},\overline{s_2}<\overline{s}=\overline{r}$. This
 means that $s_2\in S-\{s\}$ and  $r_1\in R-\{r\}$. Noting that
 $\overline{s-r}\in supp(s)\cup supp(r)$, we have either $\overline{s-r}\in supp(s)$
or $\overline{s-r}\in  supp(r)$. If $\overline{s-r}\in supp(s)$ then
$\overline{s-r}\in Irr(S-\{s\})$ which contradicts
$\overline{s-r}=c_1\overline{s_2}d_1\circ v_1$; if
$\overline{s-r}\in supp(r)$ then $\overline{s-r}\in Irr(R-\{r\})$
which contradicts $\overline{s-r}=a_1\overline{r_1}b_1\circ u_1$.
This shows that $s= r$ and then $S\subseteq R$. Similarly,
$R\subseteq S$.

Therefore, we have proved the following theorem.

\begin{theorem}\label{t5}
Let $I$ be an ideal of $kRig\langle X\rangle$ and $>$ a monomial
ordering on $Rig\langle X\rangle$. Then there exists uniquely the
reduced Gr\"{o}bner-Shirshov basis $S$ for $I$.
\end{theorem}

\section{Composition-Diamond lemma for commutative semirings}

In this section, we will give Gr\"{o}bner-Shirshov bases theory for
commutative semirings which is almost the same as the case of
semirings. The compositions of commutative semiring are simpler.

A semiring $(A, \circ, \cdot, \theta, 1)$ is commutative if  $(A,
\cdot,  1)$  is a commutative monoid. The class of  commutative
semirings is a variety.  A free commutative semiring $Rig[ X]$
generated by a set $X$ is defined as usual. Let $([X],\cdot,1)$ be
the free commutative monoid generated by $X$. If one fixes some
linear ordering $<$ on the set $[X]$, then any element of $Rig[ X]$
has an unique form $\theta$, or $w=u_1\circ u_2\circ\cdots\circ
u_n$, where $u_i\in [X], u_1\leq u_2\leq \ldots \leq u_n,\ n\geq1$.

For any $u, v\in [X]$, we denote $lcm_\cdot(u,v)$ the least common
multiple of $u$ and $v$ in $[X]$. Then there exist uniquely $a,b\in
[X]$ such that $lcm_\cdot(u,v)=au=bv$.

\begin{definition}
Let $<$ be a monomial ordering on $Rig[X]$. Let $f, g$ be two monic
polynomials in $kRig[X]$ and $\overline{f}=u_1\circ
u_2\circ\cdots\circ u_n,\  \overline{g}=v_1\circ v_2\circ\cdots\circ
v_m$ where each $u_i ,v_j\in
 [X]$. For any pair $(a, b)\in \{(a_{ij}, b_{ij})\ |\  1\leq i  \leq n, \ 1\leq j  \leq
 m\}$ where $a_{ij}, b_{ij}\in [X] $ such that $lcm_\cdot (u_i,
 v_j)=a_{ij}u_i=b_{ij}v_j$,  we call $(f,g)_{w}=af\circ u-bg\circ v$ the  composition of $f$
and $g$ with respect to $w$ where $w=lcm_\circ (a\overline{f},
b\overline{g})=a\overline{f}\circ u=b\bar{g}\circ v$.
\end{definition}

\begin{definition}
Suppose that $w$ is a monomial in  $Rig[X]$, $S$  a set of monic
polynomials in $kRig[ X]$ and $h$  a polynomial. Then  $h$ is
trivial modulo $(S,w)$, denoted by $h\equiv0 \ mod(S,w)$, if
$h=\sum\limits_{i}\alpha_ia_is_i\circ u_i$, where each
$\alpha_{i}\in{k}$, $a_i\in  [X], u_i\in Rig[X]$, $s_i\in{S}$ and
$a_i\overline{s_i}\circ u_i<w$.

The set $S$ is called a Gr\"{o}bner-Shirshov basis in $kRig[X]$ if
any composition in $S$ is trivial modulo $S$ and corresponding to
$w$.
\end{definition}

\noindent{\bf Remark}  For any given monic polynomials $f,g\in
kRig[X]$, there are finitely many compositions $(f,g)_w$. Therefore
we may use computer to realize Shirshov's algorithm to find a
Gr\"{o}bner-Shirshov basis $S^{comp}$ for a finite set $S$ in
$kRig[X]$. However, the reduced Gr\"{o}bner-Shirshov basis of
$Id(S)$ is generally infinite even if both $S$ and $X$  are finite,
see Example \ref{exa1}.

\ \

The following theorems can be similarly proved to Theorems \ref{t1}
and \ref{t5} respectively. We omit the detail.

\begin{theorem}\label{t3}{\em(Composition-Diamond lemma for commutative semirings)}\ \  Let $S$ be a set of monic
polynomials in $kRig[ X]$ and  $>$ a monomial ordering on $Rig[ X]$.
Then the following statements are equivalent.
 \begin{enumerate}
\item[(1)] $S $ is a Gr\"{o}bner-Shirshov basis in $kRig[ X]$.
\item[(2)] $ f\in Id(S)\Rightarrow \bar{f}=a\overline{s}\circ u$
for some $a\in
 [X], u\in Rig[ X]$  and $s\in S$.
\item[(2$^{'}$)]  $f\in Id(S)\Rightarrow
f=\alpha_1a_1s_1\circ u_1+\alpha_2a_2s_2\circ
u_2+\ldots+\alpha_na_ns_n\circ u_n $, where $a_1\overline{s_1}\circ
u_1>a_2\overline{s_2}\circ u_2>\ldots>a_n\overline{s_n}\circ u_n$,
$\alpha_{i}\in{k},\ a_i\in [X],\ u_i\in Rig[ X]$, $s_i\in{S}$.
\item[(3)] $Irr(S) = \{ w\in  Rig[ X] \ | \ w \neq {a\overline{s}\circ u}
 \mbox{ for  any} \ a\in
 [X], u\in Rig[ X], s\in S\}$
is a $k$-basis of $kRig[
 X|S]=kRig[
 X]/Id(S)$.
\end{enumerate}
\end{theorem}

\begin{theorem}\label{t6}
Let $I$ be an ideal of $kRig[ X]$ and $>$ a monomial ordering on
$Rig[ X]$. Then there exists uniquely the reduced
Gr\"{o}bner-Shirshov basis $S$ for $I$.
\end{theorem}

\section{Applications}

In 2004, M. Fiore and T. Leinster \cite{mar} find a strongly
normalizing reduction system and a normal form of the semiring
$\mathbb{N}[x]/(x=1+x+x^2)$. Actually $\mathbb{N}[x]/(x=1+x+
x^2)=Rig[ x|x=1\circ x\circ x^2]$. Now, we use the
Composition-Diamond lemma for commutative semirings, i.e. Theorem
\ref{t3}, to find a Gr\"{o}bner-Shirshov basis and a normal form of
this semiring.

We define a monomial  ordering on $Rig[
 x]$ first.  We order $[x]$ by degree ordering $\prec$: $x^n\prec x^m\Leftrightarrow n\prec m$.

For any $u\in Rig[
 x]$,  $u$ can be uniquely expressed as $u=u_1\circ u_2  \circ \cdots \circ u_n$,
 where $u_1, u_2,  \ldots, u_n\in [x]$, and $u_1\preceq u_2  \preceq\ldots\preceq
 u_n$.  Denote
 $$
wt(u)=(u_n, u_{n-1}, \ldots , u_1).
 $$

We order $Rig[
 x]$ as follows: for any $u,v\in Rig[x]$, if one of the sequences is not a prefix of other, then
$$
u<v\Longleftrightarrow wt(u)<wt(v)\ \mbox{ lexicographically};
$$
if the sequence of $u$ is  a prefix of the sequence of $v$, then
$u<v$.

Then, it is clear that $<$ on $Rig[
 x]$  is a monomial ordering.

\begin{theorem}\label{t2}
Let the ordering on $Rig[ x]$ be as above.  Then $kRig[ x|x=1\circ
x\circ x^2]=kRig[ x|S]$ and $S$ is a Gr\"{o}bner-Shirshov basis in
$kRig[ x]$, where  $S$ consists of the following relations
\begin{enumerate}
\item[1.]\ $x^4=1\circ1\circ x^2$,
\item[2.]\ $x\circ x^3=1\circ x^2$,
\item[3.]\ $1\circ x^2\circ x^n=x^n\ \ (1\leq n\leq 3)$.
\end{enumerate}
\end{theorem}
{\bf Proof:} We denote $i\wedge j$ the composition of the type $i$
and type $j$.

Let us check all the possible compositions.

For  $1\wedge 1$, there is no composition.

For  $1\wedge 2$, the ambiguities $w$ of all possible compositions
are: $ 1)\  x^4\circ x^6 \ \  2)\ x^2\circ x^4 $

For  $1\wedge 3$, the ambiguities $w$ of all possible compositions
are:
$$
3)\ x^4\circ x^6\circ x^{n+4}\ \  4)\ x^{2}\circ x^4\circ x^{n+2}\ \
5)\ x^{4-n}\circ x^{6-n}\circ x^4
$$
where $1\leq n\leq 3$.

For  $2\wedge 2$, the ambiguity $w$ of all possible composition is:
$6)\ x\circ x^3\circ x^5$

For  $2\wedge 3$, the ambiguities $w$ of all possible compositions
are:
$$
\begin{array}{lll}
 7) \ x\circ x^3\circ
x^5\circ x^{n+3} & 8)\    x \circ x^3\circ x^{5-n}\circ x^{3-n}  & 9)\ x^4\circ x^2\circ 1\circ x^n \\
10) \ x^{n+2}\circ x^n \circ 1\circ x^2 & 11)\ x\circ  x^3\circ
x^{n+1}&
\end{array}
$$
where $1\leq n\leq 3$.

For  $3\wedge 3$, the ambiguities $w$ of all possible compositions
are:
$$
\begin{array}{lll}
12)\ 1\circ x^2\circ x^n\circ x^{n+2}\circ x^{n+m} & 13)\ x\circ
x^3\circ x^2\circ 1\circ x^{m} & 14) \ 1\circ x^2\circ x^p\circ
x^{p-2}\circ x^{p+m-2} \\15)\ 1\circ x^2\circ x^t\circ
x^{2+t-n}\circ x^{t-n} & 16)\  1\circ x^n\circ x^2\circ x^4\circ
x^{m+2}& 17) \ 1\circ x\circ x^2\circ x^3\\ 18)\  x^4\circ x^2\circ
x^3\circ 1& 19)\ x^n\circ 1\circ x^2\circ x^m&
\end{array}
$$
where $1\leq n,m,t\leq 3$, $2\leq p \leq 3$ and $t\geq n$.

We have to check that all these compositions are trivial $mod(S,w)$.
Here, for example, we just check 1), 11), 12) and 17). Others are
similarly proved.

For 17), let $f=1\circ x\circ x^2-x,\ g=1\circ x\circ x^2-x$. Then
$w=1\circ x\circ x^2\circ x^3$  and
\begin{eqnarray*}
(f,g)_{w}&=&( 1\circ x\circ x^2-x)\circ x^3- (1\circ x\circ x^2-x)x\circ 1\\
&=&1\circ x^2-x\circ x^3\\
&\equiv&0.
\end{eqnarray*}
From this it follows that we have the relation 2.

For 11), there are three cases to consider.

Case 1. $n=1$, $f=1\circ x\circ x^2-x, \ g=x\circ x^3-1\circ x^2$.
Then $w=x\circ x^3\circ x^2$ and
\begin{eqnarray*}
(f,g)_{w}&=&( 1\circ x\circ x^2-x)x- (x\circ
x^3-1\circ x^2)\circ x^2\\
&=&1\circ x^2\circ x^2- x^2\\
&\equiv&0.
\end{eqnarray*}

Case 2.  $n=2$, $f=1\circ x^2\circ x^2-x^2, \ g=x\circ x^3-1\circ
x^2$. Then $w=x\circ x^3\circ x^3$ and
\begin{eqnarray*}
(f,g)_{w}&=&( 1\circ x^2\circ x^2-x^2)x- (x\circ
x^3-1\circ x^2)\circ x^3\\
&=&1\circ x^2\circ x^3- x^3\\
&\equiv&0.
\end{eqnarray*}
By Case 1 and Case 2 we have the relations 3.

Case 3. $n=3$, $f=1\circ x^2\circ x^3-x^3, \ g=x\circ x^3-1\circ
x^2$. Then $w=x\circ x^3\circ x^4$ and
\begin{eqnarray*}
(f,g)_{w}&=&( 1\circ x^2\circ x^3-x^3)x- (x\circ
x^3-1\circ x^2)\circ x^4\\
&=&1\circ x^2\circ x^4- x^4\\
&\equiv&1\circ x\circ x^3-x^4\\
&\equiv&1\circ 1\circ x^2-x^4\\
&\equiv&0.
\end{eqnarray*}
By Case 3 we have the relation 1.

For 12), let $f=1\circ x^2\circ x^n-x^n, \ g=1\circ x^2\circ
x^m-x^m$. Then $w=1\circ x^2\circ x^n\circ x^{n+2}\circ x^{n+m}$ and
\begin{eqnarray*}
(f,g)_{w}&=&( 1\circ x^2\circ x^n-x^n)\circ  x^{n+2}\circ x^{n+m}-
(1\circ x^2\circ
x^m-x^m)x^n\circ 1\circ x^2\\
&=&x^{n+m}\circ 1\circ x^2-x^n\circ  x^{n+2}\circ x^{n+m}\\
&\equiv&  x^{n+m}\circ 1\circ x^2-1\circ  x^{2}\circ x^{n+m}\\
&\equiv&0.
\end{eqnarray*}

For 1), let $f=x\circ x^3-1\circ x^2, \ g=x^4-1\circ1\circ x^2$.
Then $w=x^4\circ x^6$ and

\begin{eqnarray*}
(f,g)_{w}&=&(x\circ x^3-1\circ x^2)x^{3}- (x^4-1\circ1\circ x^2)\circ x^6\\
&=&1\circ1\circ x^2\circ x^6-x^3\circ x^5\\
&\equiv& 1\circ1\circ x^2\circ  x^2\circ x^2\circ x^4- 1\circ x^2 \\
&\equiv& 1\circ x^2 \circ x^2 \circ x^4- 1\circ x^2 \\
&\equiv&  x^2 \circ x^4- 1\circ x^2 \\
&\equiv&  1\circ x^2- 1\circ x^2 \\
&\equiv&0.
\end{eqnarray*}

So   $S$ is a Gr\"{o}bner-Shirshov basis in $kRig[ x]$.

The above proof implies that  $kRig[ x|x=1\circ x\circ x^2]=kRig[
x|S]$. \hfill $\blacksquare$

\ \

By
 Theorems \ref{t3} and \ref{t2}, we have  the following corollary.

\begin{corollary} (\cite{mar}) A normal form of the semiring $Rig[ x|x=1\circ x\circ
x^2]$ is the set
$$
\{1^{\circ (m+1)}\circ x^2, \ 1^{\circ m}\circ x^{\circ n}, \
1^{\circ m}\circ (x^3)^{\circ n}, \ x^{\circ m}\circ (x^2)^{\circ
n},\ (x^2)^{\circ m}\circ (x^3)^{\circ n} | n,m\geq0\}.
$$
\end{corollary}

\ \

In 1995, A. Blass \cite{blass} finds a normal form of the semiring
 $\mathbb{N}[x]/(x=1+x^2)$. Clearly, $\mathbb{N}[x]/(x=1+x^2)=Rig[ x|x=1\circ
 x^2]$. We use Theorem \ref{t3} to find a
Gr\"{o}bner-Shirshov basis  and a normal form of this semiring which
is different from \cite{blass}.

\begin{theorem}\label{t4}
Let the ordering  on $Rig[ x]$ be as in Theorem \ref{t2}. Then
$kRig[ x|x=1\circ x^2]=kRig[ x|S]$ and $S$ is a Gr\"{o}bner-Shirshov
basis in $kRig[ x]$, where $S$ consists of the following relations
\begin{enumerate}
\item[1.]\ $1\circ x^2=x$,
\item[2.]\ $x\circ x^4=1\circ x^3$,
\item[3.]\ $x^5=1\circ x^4$,
\item[4.]\ $1\circ x^3\circ x^n=x^n \ \   (3\leq n\leq 4)\}$.
\end{enumerate}
\end{theorem}
{\bf Proof:} Let us check all the possible compositions.

For  $1\wedge 1$, the ambiguity $w$ of all possible composition is:
$1)\ \ 1\circ x^2\circ x^4$

For  $1\wedge 2$, the ambiguities $w$ of all possible compositions
are:
$$
2)\ 1\circ x^2\circ x^5 \ \ \ 3)\ x^2\circ x^4\circ x \ \ \ 4)\
x^{3}\circ x\circ x^{4} \ \ \ 5)x^{6}\circ x^4\circ x
$$

For  $1\wedge 3$, the ambiguities $w$ of all possible compositions
are: $ 6)\ x^{3}\circ
 x^5\ \ \  7)\  x^7\circ x^5
$

For  $1\wedge 4$, the ambiguities $w$ of all possible compositions
are:
$$
\begin{array}{lll}
8) \ 1\circ x^2\circ x^5\circ x^{n+2} & 9)\ x \circ x^3\circ 1\circ
x^{n} & 10)\ x^{n-2} \circ x^n\circ 1\circ x^3\\  11)\ x^2\circ
1\circ x^3\circ x^{n} & 12)\ x^5\circ x^3\circ 1\circ x^n & 13)\
x^{n+2}\circ x^n \circ 1\circ x^3\\
\end{array}
$$
where $3\leq n\leq 4$.

For  $2\wedge 2$, the ambiguity $w$ of all possible composition is:
$14)\ x\circ x^4\circ x^7$

For  $2\wedge 3$, the ambiguities $w$ of all possible compositions
are: $15)\  x^2\circ x^5\ \ \  16)\ x^8\circ x^5$

For  $2\wedge 4$, the ambiguities $w$ of all possible compositions
are:
$$
\begin{array}{lll}
17)\ x\circ x^4\circ x^7\circ x^{n+4} &  18)\   x\circ x^4\circ
x^{4-n}\circ x^{7-n} & 19)\  x^6\circ x^3\circ 1\circ x^{n}\\
20)\ x^{n+3}\circ x^n\circ 1\circ x^3 & 21)\  x\circ x^4\circ
x^{n+1}&

\end{array}
$$
where $3\leq n\leq 4$.

For  $3\wedge 3$, there is no composition.

For  $3\wedge 4$, the ambiguities $w$ of all possible compositions
are:
$$
22)\ x^8\circ x^{n+5}\circ x^5 \ \ \   23)\  x^2\circ x^{n+2}\circ
x^5\ \ \  24) \  x^{5-n}\circ x^{8-n}\circ x^5
$$
where $3\leq n\leq 4$.

For  $4\wedge 4$, the ambiguities $w$ of all possible compositions
are:
$$
\begin{array}{lll}
25)\ x^{n+m}\circ x^{m+3}\circ x^m\circ x^{3}\circ 1&  26)\
x^{m-3}\circ x^{m+n-3}\circ x^m\circ x^{3}\circ 1 & 27)\ x\circ
x^{4}\circ x^4\circ x^{3}\circ 1 \\28)\ x^6\circ x^{n+3}\circ
x^3\circ 1\circ x^{m} &29)\  x^4\circ 1\circ x^3\circ x^{3}&
\end{array}
$$
where $3\leq n,m\leq 4$.

We have to check that all these compositions are trivial $mod(S,w)$.
Here, for example, we just check 1), 4), 16), 21), 23), 25) and 27).
Others are similarly proved.

For 1), let $f=1\circ x^2-x, \ g=1\circ x^2-x$. Then $w=1\circ
 x^2\circ x^4$, and
\begin{eqnarray*}
(f,g)_{w}&=&( 1\circ x^2-x)\circ x^4- (1\circ x^2-x)x^2\circ 1\\
&=&1\circ x^3-x\circ x^4\\
&\equiv&0.
\end{eqnarray*}
It follows that we have the relation 2.

For 4), let $f=1\circ x^2-x, \ g=x\circ x^4-1\circ x^3$. Then
$w=x^3\circ
 x\circ x^4$, and
\begin{eqnarray*}
(f,g)_{w}&=&( 1\circ x^2-x)x\circ x^4- (x\circ x^4-1\circ x^3)\circ x^3\\
&=&1\circ x^3\circ x^3-x^2\circ x^4\\
&\equiv&1\circ x^3\circ x^3-x^3\\
&\equiv&0.
\end{eqnarray*}
Then we have the relation 4 for $n=3$.

For 21), there are two cases to consider.

Case 1. $n=3$,  let $f=x\circ x^4-1\circ x^3, \ g=1\circ x^3\circ
x^3- x^3$. Then $w=x\circ x^4\circ
 x^4$, and
\begin{eqnarray*}
(f,g)_{w}&=&( x\circ x^4-1\circ x^3)\circ x^4- (1\circ
x^3\circ x^3- x^3)x\\
&=& x^4-1\circ x^3\circ x^4\\
&\equiv&0.
\end{eqnarray*}
Then we have the relation 4 for $n=4$.

Case 2. $n=4$,  let $f=x\circ x^4-1\circ x^3, \ g=1\circ x^3\circ
x^- x^4$. Then $w=x\circ x^4\circ
 x^5$, and
\begin{eqnarray*}
(f,g)_{w}&=&( x\circ x^4-1\circ x^3)\circ x^4- (1\circ
x^3\circ x^4- x^4)x\\
&=& x^5-1\circ x^3\circ x^5\\
&\equiv& x^5-1\circ x^4\\
&\equiv&0.
\end{eqnarray*}
Then we have the relation 3.

For 16), let $f=x\circ x^4-1\circ x^3, \ g= x^5-1\circ x^4$. Then
$w=x^8\circ x^5$, and
\begin{eqnarray*}
(f,g)_{w}&=&(x\circ x^4-1\circ x^3) x^4- (x^5-1\circ x^4)\circ x^8\\
&=&1\circ x^4\circ x^8- x^4\circ x^7\\
&\equiv&1\circ x^4\circ x^3\circ x^7- x^4\circ x^7\\
&\equiv&x^4\circ x^7- x^4\circ x^7\\
&\equiv&0.
\end{eqnarray*}

For 23), let $f=x^5- 1\circ x^4, \ g=1\circ x^3\circ x^n- x^n $.
Then $x^2\circ x^{n+2}\circ x^5$, and
\begin{eqnarray*}
(f,g)_{w}&=&(x^5- 1\circ x^4) \circ x^2\circ x^{n+2}- (1\circ x^3\circ x^n- x^n)x^2\\
&=& x^{n+2}- 1\circ x^4 \circ x^2\circ x^{n+2} \\
&\equiv&x^{n+2}- 1\circ  x^3\circ x^{n+2} \\
\end{eqnarray*}
There are two cases to consider.

Case 1. $n=3$. We have
\begin{eqnarray*}
(f,g)_{w}&\equiv&x^{n+2}- 1\circ  x^3\circ x^{n+2}\\
&\equiv&x^{5}- 1\circ x^3\circ x^{5}\\
&\equiv&1\circ x^{4}- 1\circ x^3\circ 1\circ x^{4}\\
&\equiv&1\circ x^{4}- 1\circ x^{4}\\
&\equiv&0.
\end{eqnarray*}

Case 2. $n=4$. We have
\begin{eqnarray*}
(f,g)_{w}&\equiv&x^{n+2}- 1\circ  x^3\circ x^{n+2}\\
&\equiv&x^{6}- 1\circ x^3\circ x^{6}\\
&\equiv& x\circ x^{5}- 1\circ x^3\circ x\circ x^{5}\\
&\equiv& x\circ x^{5}- 1\circ x^2\circ x^{5}\\
&\equiv& x\circ x^{5}- x\circ x^{5}\\
&\equiv&0.
\end{eqnarray*}

For 25), let $f=1\circ x^3\circ x^n- x^n, \ g= 1\circ x^3\circ x^m-
x^m$. Then $w=x^{n+m}\circ x^{m+3}\circ x^m\circ x^{3}\circ 1 $, and
\begin{eqnarray*}
(f,g)_{w}&=&(1\circ x^3\circ x^n- x^n) x^m\circ 1\circ x^3- (1\circ
x^3\circ x^m-
x^m)\circ x^{n+m}\circ x^{3+m}\\
&=&x^m\circ x^{n+m}\circ x^{3+m}- x^{n+m}\circ 1\circ x^3\\
&\equiv&x^mx^n- 1\circ x^3\circ x^{n+m}\\
&\equiv&x^{n+m}- 1\circ x^3\circ x^{n+m}.\\
\end{eqnarray*}
There are three cases to consider.

Case 1. $n=m=3$. We have

\begin{eqnarray*}
(f,g)_{w}&\equiv&x^mx^n- 1\circ x^3\circ x^{n+m}\\
&\equiv&x^{6}- 1\circ x^3\circ x^{6}\\
&\equiv& x\circ x^{5}- 1\circ x^3\circ x\circ x^{5}\\
&\equiv& x\circ x^{5}- 1\circ x^2\circ x^{5}\\
&\equiv& x\circ x^{5}- x\circ x^{5}\\
&\equiv&0.
\end{eqnarray*}

Case 2.  $n=m=4$. We have

\begin{eqnarray*}
(f,g)_{w}&\equiv&x^mx^n- 1\circ x^3\circ x^{n+m}\\
&\equiv&x^{8}- 1\circ x^3\circ x^{8}\\
&\equiv& x^3\circ x^{7}- 1\circ x^3\circ x^3\circ x^{7}\\
&\equiv& x^3\circ x^{7}-  x^3\circ x^{7}\\
&\equiv&0.
\end{eqnarray*}

Case 3. $n=3, m=4$, or $m=3, n=4$. We have

\begin{eqnarray*}
(f,g)_{w}&\equiv&x^mx^n- 1\circ x^3\circ x^{n+m}\\
&\equiv&x^{7}- 1\circ x^3\circ x^{7}\\
&\equiv& x^2\circ x^{6}- 1\circ x^3\circ x^2\circ x^{6}\\
&\equiv& x^2\circ x^{6}-  x^3\circ x\circ x^{6}\\
&\equiv& x^2\circ x^{6}-  x^2\circ x^{6}\\
&\equiv&0.
\end{eqnarray*}

For 27), let $f=1\circ x^3\circ x^3- x^3, \ g= 1\circ x^3\circ x^4-
x^4$. Then $w=x\circ x^{4}\circ x^4\circ x^{3}\circ 1$, and
\begin{eqnarray*}
(f,g)_{w}&=&(1\circ x^3\circ x^3- x^3) x\circ 1\circ x^3- (1\circ
x^3\circ x^4-
x^4)\circ x\circ x^4\\
&=&x^4\circ x\circ x^4- x^4\circ 1\circ x^3 \\
&\equiv&1\circ x^3\circ x^4- x^4\circ 1\circ x^3 \\
&\equiv&0.
\end{eqnarray*}

Therefore  $S$ is a Gr\"{o}bner-Shirshov basis in $kRig[ x]$.

The above proof implies that  $kRig[ x|x=1\circ x^2]=kRig[ x|S]$.

We complete the proof. \hfill $\blacksquare$

\ \

By Theorems \ref{t3} and \ref{t4}, we have  the following corollary.

\begin{corollary}\label{coro5.4}  A normal form of the semiring $Rig[ x|x=1\circ
x^2]$ is the set
$$
\{ (1^{\circ n} \circ x^{\circ m})x^t,  1^{\circ n}\circ x^3,
1^{\circ n}\circ (x^4)^{\circ m}\ |\ n,m\geq 0, \ 0\leq t\leq 3\}.
$$
\end{corollary}

In order to compare an another normal form of the semiring $Rig[
x|x=1\circ x^2]$
$$
\{ 1^{\circ n} \circ x^2\circ x^4, \ 1^{\circ n} \circ (x^2)^{\circ
m}, \ (x^2)^{\circ m}\circ (x^4)^{\circ t},\ 1^{\circ n} \circ
(x^4)^{\circ t}\mid \ n,m,t\geq 0\}
$$
given by A. Blass \cite{blass}, we need the following lemma.

\begin{lemma}\label{l2} Suppose that $\Gamma,\Sigma$ are two subsets
of the semiring $Rig\langle X\rangle$ and $\rho$ is a congruence on
$Rig\langle X \rangle$. Suppose that $\Gamma$ is a normal form of
$Rig\langle X|\rho \rangle$. If $f:\Gamma\rightarrow \Sigma$ is a
bijective mapping such that for any $u\in \Gamma$, $f(u)\rho=u\rho$,
then $\Sigma$ is also a normal form of $Rig\langle X|\rho\rangle$.
\end{lemma}
{\bf Proof:} For any $u\in Rig\langle X \rangle$, since $\Gamma$ is
a normal form of the semiring $Rig\langle X|\rho\rangle$, there is
uniquely $v\in \Gamma$, such that $u\rho=v\rho$. Hence
$u\rho=f(v)\rho$, where $f(v)\in \Sigma$.

For any two different $u,v\in\Sigma$, if $u\rho=v\rho$, then
$f^{-1}(u)\rho=f^{-1}(v)\rho$ and hence $f^{-1}(u)\neq f^{-1}(v)$, a
contradiction. This shows that $\Sigma$ is a normal form of the
semiring $Rig\langle X|\rho\rangle$. \hfill $\blacksquare$

\ \

\begin{corollary}(\cite{blass})
A normal form of the semiring $Rig[ x|x=1\circ x^2]$ is the set
$$
\{ 1^{\circ n} \circ x^2\circ x^4, \ 1^{\circ n} \circ (x^2)^{\circ
m}, \ (x^2)^{\circ m}\circ (x^4)^{\circ t},\ 1^{\circ n} \circ
(x^4)^{\circ t}\mid \ n,m,t\geq 0\}.
$$
\end{corollary} {\bf Proof:} We denote
$$
\Gamma=\{ (1^{\circ n} \circ x^{\circ m})x^t,  1^{\circ n}\circ x^3,
1^{\circ n}\circ (x^4)^{\circ m}\ |\ n,m\geq 0, \ 0\leq t\leq 3\},
$$
$$
\Sigma=\{ 1^{\circ n} \circ x^2\circ x^4, \ 1^{\circ n} \circ
(x^2)^{\circ m}, \ (x^2)^{\circ m}\circ (x^4)^{\circ t},\ 1^{\circ
n} \circ (x^4)^{\circ t}\mid \ n,m,t\geq 0\}
$$
and $\rho$ the congruence on $Rig[ x]$ generated by $\{x=1\circ
x^2\}$.

Define
\begin{eqnarray*}
f:&&\Gamma\rightarrow \Sigma,\\
&&1^{\circ n} \circ x^{\circ m}\mapsto 1^{\circ (n+m)} \circ
(x^2)^{\circ m},\\
&&(1^{\circ n} \circ x^{\circ m})x\mapsto 1^{\circ n} \circ
(x^2)^{\circ (n+m)},\\
&&(1^{\circ n} \circ x^{\circ m})x^2\mapsto (x^2)^{\circ (n+m)}
\circ
(x^4)^{\circ m},\\
&&(1^{\circ n} \circ x^{\circ m})x^3\mapsto (x^2)^{\circ n}
 \circ(x^4)^{\circ (n+m)},\\
&&1^{\circ n}\circ x^3\mapsto 1^{\circ n}\circ x^2\circ x^4,\\
&&1^{\circ n}\circ (x^4)^{\circ m}\mapsto 1^{\circ n}\circ
(x^4)^{\circ m}.
\end{eqnarray*}
Then  $f$ is a bijective mapping and for any $u\in \Gamma$,
$f(u)\rho=u\rho$ since $f(u)$ is obtained by $u$ replacing $x,\ x^3$
for $1\circ x^2,\ x^2\circ x^4$ respectively.

Now the result follows from Corollary \ref{coro5.4} and Lemma
\ref{l2}.
 \hfill $\blacksquare$

\ \

Let $X$ be a well ordered set, $Z$ the integer ring and $Z\langle
X\rangle$ the semigroup ring over $Z$. It is easy to see that
$(Z\langle X\rangle, \circ, \cdot) $ is a semiring with the
operations  $ f\circ g:=f+g, f\cdot g:=f\times g$, where $f, g$ are
polynomials in $Z\langle X\rangle$. Now, we represent the semiring
$Z\langle X\rangle$ by generators and defining relations.

Let $X^{-1}=\{x^{-1}| x\in X\}$. We define a monomial ordering on
$Rig\langle X\cup X^{-1}\cup 1^{-1}\rangle$ first.

For any $x, y \in X$, we define $x^{-1}>x>y^{-1}>y$ if $x>y$ and
$x>1^{-1}>1$. Then we define the inverse deg-lex ordering $\preceq$
on $ \{X\cup X^{-1}\cup 1^{-1}\}^*$.

For any $u\in Rig\langle X\cup X^{-1}\cup 1^{-1}\rangle$,  $u$ can
be uniquely expressed as $u=u_1\circ u_2  \circ \cdots \circ u_n$,
where $u_1, u_2,  \ldots, u_n\in \{X\cup X^{-1}\cup 1^{-1}\}^*$ and
$u_1\preceq u_2  \preceq\ldots\preceq
 u_n$.  Denote
 $$
wt(u)=(deg(u), n, u_1, u_{2}, \ldots , u_n).
 $$

We order $Rig\langle
 X\cup X^{-1}\cup 1^{-1}\rangle$ as follows: for any $u,v\in Rig\langle X\cup X^{-1}\cup 1^{-1}\rangle$,

$$
u>v\Longleftrightarrow wt(u)>wt(v)\ \mbox{ lexicographically}.
$$

Then, it is clear that $>$ on $Rig\langle X\cup X^{-1}\cup
1^{-1}\rangle$  is a monomial ordering.

\begin{theorem}
Let the ordering be as above. Then $Z\langle X\rangle\cong
Rig\langle
 X\cup X^{-1}\cup 1^{-1}| S\rangle$ as semirings and a Gr\"{o}bner-Shirshov basis $S$
in $kRig\langle
 X\cup X^{-1}\cup 1^{-1}\rangle$  consists of the following
 relations:
\begin{enumerate}
\item[1.]\ $x\circ x^{-1}=\theta$,
\item[2.]\ $1\circ 1^{-1}=\theta$,
\item[3.]\ $ x^{-1}y^{-1}=xy$,
\item[4.]\ $xy^{-1}=x^{-1}y$,
\item[5.]\ $x^\epsilon 1^{-1}
=x^{-\epsilon}$,
\item[6.]\ $1^{-1}x^\epsilon
=x^{-\epsilon}$,
\end{enumerate}
where $x, y\in X,  \ \epsilon=\pm 1$. As a result, a normal form of
the
 semiring $Rig\langle
 X\cup X^{-1}\cup 1^{-1}\ | S \rangle$ is the set
\begin{eqnarray*}
Irr(S)&=&\{x_{11}^{\epsilon_1}x_{12}\cdots x_{1n_1}\circ
x_{21}^{\epsilon_2}x_{22}\cdots
 x_{2n_2} \circ \cdots \circ x_{m1}^{\epsilon_m}x_{m2}\cdots
 x_{mn_m}\\
 &&\ \ \ \ | \ x_{ij}\in X, m\geq 0, \epsilon_i=\pm 1, n_i\geq 0, i=1, \ldots, m
 \}
 ,
\end{eqnarray*}
where $x_{i1}^{\epsilon_1}x_{i2}\cdots
x_{in_i}=1^{\epsilon_1}$ if $n_i=0$.
\end{theorem}
{\bf Proof:} It is easy to see that
$$
\sigma: Z\langle X\rangle\rightarrow Rig\langle
 X\cup X^{-1}\cup 1^{-1}| S\rangle, \ \ \epsilon x_{i1}x_{i2}\cdots x_{it}\mapsto
 x_{i1}^{\epsilon}x_{i2}\cdots x_{it}, \ 0\mapsto\theta
$$
is a semiring isomorphism, where $\epsilon=\pm1$. Since
$Irr(S)=\sigma(Z\langle X\rangle)$, $Irr(S)$ is a $k$-basis of
$kRig\langle
 X\cup X^{-1}\cup 1^{-1}\ | S \rangle$. Therefore, by using
Theorem \ref{t1}, $S$ is a Gr\"{o}bner-Shirshov basis in
$kRig\langle
 X\cup X^{-1}\cup 1^{-1}\rangle$. \hfill $\blacksquare$

\ \

Let $(\mathbb{N}, \circ, \cdot)$ be the natural numbers semiring,
where for any $n,m\in \mathbb{N}$, $n\circ m:=n+m,\ n\cdot
m:=n\times m$. Then $(\mathbb{N}, \circ, \cdot)=Rig[ x \ |\ x=1]$.
For any congruence $\rho $ on $\mathbb{N}$, we have
$\mathbb{N}/\rho=Rig[ x \ |\  x=1, \rho]$. Let the ordering on $Rig[
x]$ be defined as in Theorem \ref{t2}. By Shirshov algorithm, we are
able to find a Gr\"{o}bner-Shirshov basis $\{x=1, \rho\}^{comp}$ for
the set $\{x=1, \rho\}$. Suppose $\{x=1\}\cup S=\{x=1,
\rho\}^{comp}$. Then by Theorem \ref{t6}, we may assume that
$\{x=1\}\cup S$ is the reduced Gr\"{o}bner-Shirshov basis. Since
$\{x=1\}\cup S$ is minimal,   each element in $S$ has the form
$1^{\circ n}=1^{\circ m}, n>m\in \mathbb{N}$ and  $S$ contains only
one element, say, $1^{\circ n}=1^{\circ m}, n>m\in \mathbb{N}$. It
follows that the congruence $\rho$ on $\mathbb{N}$ is generated by
one element $(n,m)$.

Thus, we have the following corollary.

\begin{corollary}
Each congruence on the semiring $\mathbb{N}$  is generated by one
element. In particular, $\mathbb{N}$ is Noetherian.
\end{corollary}

\ \

For a  commutative algebra $k[X|S]$ with $|X|<\infty$, it is well
known that a reduced Gr\"{o}bner-Shirshov basis of $k[X|S]$ must be
finite. It is also well known that if the ring $R$ is Noetherian
then so is the polynomial ring $R[X]$ if $|X|<\infty$. However, it
is not the case for the semiring $\mathbb{N}[x]$.

\begin{example}\label{exa1}
Considering the semiring $\mathbb{N}[x]/(x+1=x)=Rig[x | \ x\circ
1=x]$, it is easy to have that $kRig[ x | \ x\circ 1=x]=kRig[x | S
]$ where $S=\{x^n\circ1=x^n | \ n\geq 1\}$  is the reduced
Gr\"{o}bner-Shirshov basis in $kRig[
 x]$ with the ordering in Theorem \ref{t2}.
\end{example}

Now, we construct an ascending  chain of ideals in $kRig[ x ]$ as
follows.
$$
I_1 \subseteq I_2 \subseteq \ldots \subseteq I_n\subseteq \ldots
$$
where $I_n=Id(x\circ1, x^2\circ1,\ldots,  x^n\circ1)$.

For any $n\geq 1$, $x^{n+1}\circ 1\not\in I_n$. Otherwise, there
exist  $ n\geq i\geq 1, a,b\in [x], u\in Rig[ x]$ such that
$x^{n+1}\circ 1=a(x^i\circ 1)b\circ u$. This is a contradiction
because $S$ is a minimal Gr\"{o}bner-Shirshov basis in $kRig[
 x]$. Hence
$$
I_1 \subsetneq I_2  \subsetneq \ldots \subsetneq I_n\subsetneq
\ldots.
$$

Let us define  congruence relation $\rho_n$  on $\mathbb{N}[x]$
generated by the set
$$
\{(x^i\circ1, x\circ1),\ 2\leq i\leq n\}.
$$

Since $(x^{n+1}\circ1,x\circ1)\not\in \rho_{n}$, we have an infinite
ascending chain of congruences
$$
\rho_1 \subsetneq \rho_2  \subsetneq \ldots \subsetneq
\rho_n\subsetneq \ldots.
$$
Thus, we have  the following corollary.
\begin{corollary}
$\mathbb{N}[ x]$ is not Noetherian.
\end{corollary}

\ \

\noindent{\bf Acknowledgement}: We are grateful to Marcelo Fiore who
took our attention to his and T. Leinster's paper \cite{mar}.

\end{document}